\documentclass{article}
\usepackage{amsfonts}
\newtheorem{theorem}{Theorem}

\newtheorem{definition}[theorem]{Definition}

\newtheorem{lemma}[theorem]{Lemma}
\newtheorem{proposition}[theorem]{Proposition}

\newenvironment{proof}[1][Proof]{\textbf{#1.} }{\ \rule{0.5em}{0.5em}}

\title{Hausdorff Continuous Solutions of Nonlinear PDEs through the Order Completion Method}
\author{Roumen Anguelov and Elemer E Rosinger\\
Department of Mathematics and Applied Mathematics\\
University of Pretoria\\
SOUTH AFRICA\\
anguelov@scientia.up.ac.za\\
rosinger@scientia.up.ac.za}
\date{}

\begin{document}
\maketitle
\begin{abstract}
It was shown in \cite{Rosinger} that very large classes of
nonlinear PDEs have solutions which can be assimilated with usual
measurable functions on the Euclidean domains of definition of the
respective equations. In this paper the regularity of these
solutions has significantly been improved by showing that they can
in fact be assimilated with Hausdorff continuous functions. The
method of solution of PDEs is based on the Dedekind order
completion of spaces of smooth functions which are defined on the
domains of the given equations.
\end{abstract}

\section{Introduction}

The following significant {\it threefold} breakthrough was
obtained in \cite{Rosinger} with respect to solving large classes
of nonlinear PDEs, see MR 95k:35002. Namely:

\begin{itemize}
\item[a)] arbitrary nonlinear PDEs of the form
\begin{equation}\label{PDE}
T(x,D)u=f(x),{\rm  \ \ }x\in\Omega,%
\end{equation}
where
\begin{equation}\label{PDE2}
T(x,D)u=g(x,u(x),~.~.~.~,D_{x}^{p}u(x),~.~.~.~),~~
p\in\mathbb{N}^{n},\ |p|\leq m,
\end{equation}
with $g$ jointly continuous in all it arguments, $f$ in a class of
measurable functions, $\Omega \subseteq {\bf R}^n$ arbitrary open,
$m \in \mathbb{N}$ arbitrary given, and the unknown function $u :
\Omega \longrightarrow \mathbb{R}$, were proven to have%

\item[b)] solutions $u$ which can be assimilated with usual
measurable functions on $\Omega$, and

\item[c)] the solution method was based on the Dedekind order
completion of suitable spaces of smooth functions on $\Omega$.

\end{itemize}

In fact, the conditions at a)~ can further be relaxed by assuming
that $g$ may admit certain {\it discontinuities}, namely, that it
is continuous only on $(\Omega\setminus\Sigma )\times{\bf
R}^{m^*}$, where $\Sigma$ is a closed, nowhere dense subset of
$\Omega$, while $m^*$ is the number of arguments in $g$ minus $n$.
This relaxation on the continuity of $g$ may be significant since
such subsets of discontinuity $\Sigma$ can have arbitrary large
positive Lebesgue measure.

In this way, the solutions of the unprecedented large class of
nonlinear PDEs in (1) can be obtained {\it without} the use of any
sort of distributions, hyperfunctions, generalized functions, or
of methods of functional analysis. Moreover, one obtains a
general, {\it blanket regularity}, given by the fact that the
solutions constructed can be assimilated with usual measurable
functions on the corresponding domains $\Omega$ in Euclidean
spaces.

In this paper we discuss a further significant improvement of the
above mentioned results with respect to the {\it regularity}
properties of the solutions. Namely, this time we show that they
can be assimilated with the significantly smaller class of
Hausdorff continuous functions on the open domains $\Omega$. This
improvement follows, among others, from a recent breakthrough, see
\cite{QM}, which solves a long outstanding problem related to the
Dedekind order completion of spaces $C(X)$ of real valued
continuous functions on rather arbitrary topological spaces $X$.

The Hausdorff continuous functions are not unlike the usual real
valued continuous functions. For instance, they assume real values
on a dense subset of the domain and are completely determined by
the values on this subset. However, these functions may also
assume interval values on a certain subset of the domain. Hence
the concept of Hausdorff continuity is formulated within the realm
of interval valued functions. We denote by $\mathbb{A}(\Omega )$
the set of all functions defined on an open set
$\Omega\subset\mathbb{R}^{n}$ with values which are finite or
infinite closed real intervals, that is,
\[
\mathbb{A}(\Omega )=\{f:\Omega \rightarrow
\mathbb{I\,\overline{\mathbb{R}} }\},
\]%
where
$\mathbb{I\,}\overline{\mathbb{R}}=\{[\underline{a},\overline
{a}]:\underline{a},\overline{a}\in\overline{\mathbb{R}}=\mathbb{R\cup
\{\pm\infty\}},$ $\underline{a}\leq\overline{a}\}$. Given an
interval
$a=[\underline{a},\overline{a}]\in\mathbb{I\,\overline{\mathbb{R}}}$,
\[
w(a)=\left\{
\begin{tabular}
[c]{lll}%
$\overline{a}-\underline{a}$ & if & $\underline{a},\overline{a}$ finite,\\
$+\infty$ & if & $\underline{a}<\overline{a}=+\infty$ or $\underline{a}%
=-\infty<\overline{a}$,\\
0 & if & $\underline{a}=\overline{a}=\pm\infty,$%
\end{tabular}
\right.
\]
is the width of $a$, while
$|a|=\max\{|\underline{a}|,|\overline{a}|\}$ is the modulus of
$a$. An extended real interval $a$ is called proper if $w(a)>0$
and degenerate or point if $w(a)=0$. Identifying $a\in
\overline{\mathbb{R}}$ with the degenerate interval $[a,a]\in
\mathbb{I}\,\overline{\mathbb{R}}$, we consider
$\overline{\mathbb{R}}$ as a subset of
$\mathbb{I}\,\overline{\mathbb{R}}$. In this way
$\mathbb{A}(\Omega )$ contains the set of extended real valued functions, namely,%
\[
\mathcal{A}(\Omega)=\{f:\Omega \rightarrow
\overline{\mathbb{R}}\}.
\]
A partial order which extends the total order on
$\overline{\mathbb{R}}$ can be defined on
$\mathbb{I\,}\overline{\mathbb{R}}$ in more than one way. However,
it will prove useful to consider on
$\mathbb{I\,}\overline{\mathbb{R}}$ the partial order $\leq$ defined by%
\begin{equation}
\lbrack\underline{a},\overline{a}]\leq\lbrack\underline{b},\overline
{b}]\Longleftrightarrow\underline{a}\leq\underline{b},\;\overline{a}%
\leq\overline{b}. \label{iorder}%
\end{equation}
The partial order induced on $\mathbb{A}(\Omega)$ by
(\ref{iorder}) in a point-wise
way, i.e.,%
\begin{equation}
f\leq g\Longleftrightarrow f(x)\leq g(x),\;x\in \Omega, \label{forder}%
\end{equation}
is an extension of the usual point-wise order in the set of
extended real valued functions $\mathcal{A}(\Omega)$.

The application of Hausdorff continuous functions to problems in
Analysis, e.g. \cite{QM}, and to nonlinear PDEs as in this paper,
are based on the quite extraordinary fact that the set
$\mathbb{H}(\Omega)$ of all Hausdorff continuous functions on
$\Omega$ is order complete while some of its important subsets are
Dedekind order complete. We can recall that the usual spaces of
real valued functions considered in Analysis or Functional
Analysis, e.g. spaces of continuous functions, Lebesgue spaces,
Sobolev spaces, are with very few exceptions neither order
complete nor Dedekind order complete.

The definition of the concept of Hausdorff continuity and related
terminology are discussed in Section 2. The Baire operators and
the graph completion operator which are instrumental for the
definition and the properties of Hausdorff continuous functions
are also discussed in that section. In order to improve the
readability of the paper, a short account of some basic properties
of the Hausdorff continuous functions is given in the Appendix.

The use of extended real intervals in the definition of the set
$\mathbb{A}(\Omega)$ is partially motivated by the fact that the
Baire operators involve infimums and supremums which might not
exists in the realm of the usual (finite) real intervals. However,
the main motivation with regard to the present exposition is the
need to accommodate solutions of PDEs which are discontinuous at
certain points of the domain $\Omega$ and unbounded in the
neighborhood of these points, e.g. the so called finite time
blow-up. For this purpose it will prove sufficient to consider
only the nearly finite functions.
\begin{definition}\label{Defnearlyfinite}
A function $f\in\mathbb{A}(\Omega)$ is called nearly finite if
there exists an open and dense subset $D$ of $\Omega$ such that
\[
|f(x)|<+\infty,\;x\in D
\]
\end{definition}
After a brief introduction to the order completion method in
Section 3 we give the main result of the paper in Section 4,
namely, that the solutions of the equation (\ref{PDE}) through the
order completion method can be assimilated with nearly finite
Hausdorff continuous functions.

\section{Baire operators, graph completion operator and Hausdorff continuity}

For every $x\in \Omega$, $B_{\delta }(x)$ denotes the open $\delta $%
-neighborhood of $x$ in $\Omega $, that is,
\[
B_{\delta }(x)=\{y\in \Omega :||x-y||<\delta \}.
\]%
Let $D$ be a dense subset of $\Omega$.  The pair of mappings
$I(D,\Omega,\cdot),$ $S(D,\Omega,\cdot):\mathbb{A}(D )\rightarrow
\mathcal{A}(\Omega )$ defined by
\begin{eqnarray}
I(D,\Omega ,f)(x) &=&\sup_{\delta >0}\inf \{z\in f(y):y\in
B_{\delta
}(x)\cap D\},x\in \Omega ,  \label{lbfgen} \\
S(D,\Omega ,f)(x) &=&\inf_{\delta >0}\sup \{z\in f(y):y\in
B_{\delta }(x)\cap D\},x\in \Omega ,  \label{ubfgen}
\end{eqnarray}%
are called lower Baire and upper Baire operators, respectively.
Clearly for every $f\in \mathbb{A}(D )$ we have
\[
I(D,\Omega ,f)(x)\leq f(x)\leq S(D,\Omega ,f)(x),\;x\in\Omega.
\]
Hence the mapping $F:\mathbb{A}(D)\rightarrow \mathbb{A}(\Omega
)$, called a graph completion operator, where
\begin{equation}
F(D,\Omega ,f)(x)=[I(D,\Omega ,f)(x),S(D,\Omega ,f)(x)],\;x\in
\Omega ,\;f\in \mathbb{A}(\Omega ), \label{gcgen}
\end{equation}
is well defined and we have the inclusion
\begin{equation}\label{fFincl}
f(x)\subseteq F(f)(x),\ x\in\Omega.
\end{equation}
The name of this operator is derived from the fact that
considering the graphs of $f$ and $F(D,\Omega,f)$ as subsets of
the topological space $\Omega\times \overline{\mathbb{R}}$, the
graph of $F(D,\Omega,f)$ is the minimal closed set which is a
graph of interval function on $\Omega$ and contains the the graph
of $f$. In the case when $D=\Omega$ the sets $D$ and $\Omega$ will
be usually omitted from the operators' argument lists, that is,
\[
I(f)=I(\Omega,\Omega,f)\ ,\ \ S(f)=S(\Omega,\Omega,f)\ ,\ \
F(f)=F(\Omega,\Omega,f)
\]

Let us note that, the graph completion operator is monotone about
inclusion with respect to the functional argument, that is, if
$f,g\in \mathbb{A}(D)$ where $D$ is dense in $\Omega$ then
\begin{equation}\label{Finclmon1}
f(x)\subseteq g(x),\;x\in D\Longrightarrow F(D,\Omega
,f)(x)\subseteq F(D,\Omega ,g)(x),\;x\in \Omega .
\end{equation}
Furthermore, the graph completion operator is monotone about
inclusion with respect to the set $D$ in the sense that if $D_{1}$
and $D_{2}$ are dense subsets of $\Omega $ and $f\in
\mathbb{A}(D_{1}\cup D_{2})$ then
\begin{equation}\label{Finclmon2}
D_{1}\subseteq D_{2}\Longrightarrow F(D_{1},\Omega ,f)(x)\subseteq
F(D_{2},\Omega ,f)(x),x\in \Omega .
\end{equation}
This, in particular, means that for any dense subset $D$ of
$\Omega $ and $f\in
\mathbb{A}(\Omega )$ we have%
\begin{equation}\label{Fincmon3}
F(D,\Omega ,f)(x)\subseteq F(f)(x),x\in \Omega .
\end{equation}

Let $f\in \mathbb{A}(\Omega).$ For every $x\in \Omega$ the value
of $f$ is an interval $[\underline{f}(x),\overline{f}(x)]$. Hence,
the function $f$ can be written
in the form $f=[\underline{f},\overline{f}]$ where $\underline{f},\overline{f%
}\in \mathcal{A}(X)$ and $\underline{f}\leq \overline{f}$. The
lower and upper Baire operators as well as the graph completion
operator of an interval valued function
$f=[\underline{f},\overline{f}]\in \mathcal{A}(\Omega)$
can be conveniently represented in terms of the functions $\underline{f}$ and $%
\overline{f}$:
\[
I(D,\Omega,f)=I(D,\Omega,\underline{f})\;,\;\;
S(D,\Omega,f)=S(D,\Omega,\overline{f}),
\]
\[F(D,\Omega,f)=[I(D,\Omega,\underline{f}),S(D,\Omega,\overline{f})].
\]

\begin{definition}\label{DefHcont}%
A function $f\in\mathbb{A}(\Omega)$ is called Hausdorff
continuous, or H-continuous, if for every function
$g\in\mathbb{A}(\Omega)$ which satisfies the inclusion
$g(x)\subseteq f(x),$ $x\in \Omega$, we have $F(g)(x)=f(x),$ $x\in
\Omega $.
\end{definition}

The concepts of Hausdorff continuity is strongly connected to the
concepts of semi-continuity of real functions. We have the
following characterization of the fixed points of the lower and
the upper Baire operators, see \cite{Baire}:
\begin{eqnarray}
I(f) &  =f\Longleftrightarrow f{\rm \;-\; lower\;
semi\!-\!continuous\; on\;
}\Omega\label{flsc}\\
S(f) &  =f\Longleftrightarrow f{\rm \;-\; upper\;
semi\!-\!continuous\; on\;
}\Omega\label{fusc}%
\end{eqnarray}
Hence an interval function $f=[\underline{f},\overline{f}]$ is
H-continuous if and only if the following three conditions hold
\begin{itemize}
\item[(i)] $\underline{f}$ is lower semi-continuous%
\item[(ii)] $\overline{f}$ is upper semi-continuous%
\item[(iii)] the set $\{\phi\in {\cal
A}(\Omega):\underline{f}\leq\phi\leq\overline{f}\}$ does not
contain lower or upper semi-continuous functions other than
$\underline{f}$ and $\overline{f}$
\end{itemize}

The concept of H-continuity can be considered as a generalization
of the concept of continuity of real functions in the sense that
the only real (point valued) functions contained in
$\mathbb{H}(\Omega)$ are the continuous functions, that is,
\begin{equation}
\left.\begin{tabular}{c} $f\in \mathcal{A}(\Omega)$\\ $f$ is
H-continuous
\end{tabular}\right\}\Longrightarrow f\ {\rm\ is\ continuous}
\end{equation}
The H-continuous functions retain some essential properties of the
usual real continuous functions as stated in Theorem \ref{tindent}
in the Appendix. Further links with the real continuous functions
are presented in Theorems \ref{tcont} and \ref{textend}. We should
also note that any Hausdorff continuous function $f$ is
"essentially" point valued in the sense that it assumes point
values everywhere except on a set $W_f$ which is of first Baire
category, see Theorem \ref{FBC}. Through an application of the
Baire category theorem this implies that the complement of $W_f$
in $\Omega$ is a set of second Baire category. Hence %
\begin{equation}\label{Dfdense}
D_f=\Omega\setminus
W_f=\{x\in\Omega:f(x)\in\overline{\mathbb{R}}\}\ \ {\rm is\ dense\
in\ } \Omega.
\end{equation}

\section{Order completion method for nonlinear PDEs.}

The order completion method in solving general nonlinear systems
of PDEs of the form (\ref{PDE}) is based on certain very simple,
even if less than usual, approximation properties, see
\cite{Rosinger}. To give an idea about the ways the order
completion method works, we mention some of these approximations
here.

The differential operator $T(x, D)$ o n the left hand side of
(\ref{PDE}) has the following basic approximation property :

\begin{lemma}
\[ \begin{array}{l}
                 \forall~~ x_0 \in \Omega,~~ \epsilon > 0 ~~: \\
                 \\
                 \exists~~ \delta > 0,~~ P ~~\mbox{polynomial in}~~ x \in
                                                                {\bf R}^n~~ :
                                                                \\ \\
                 ~~~ | | x - x_0 | | ~\leq~ \delta ~~\Longrightarrow~~ f(x) - \epsilon ~\leq~
                                                T(x,D) P(x) ~\leq~ f(x)
   \end{array}\]
\hfill $\Box$
\end{lemma}

Consequently, we obtain :

\begin{proposition}\label{prop}
\[~~~ \begin{array}{l}
              \forall~~ \epsilon > 0~~ : \\ \\
              \exists~~ \Gamma_\epsilon \subset \Omega ~~\mbox{closed,~ nowhere~
                  dense~ in}~~  \Omega,~~ U_\epsilon \in C^\infty (\Omega)~~ :
                  \\ \\
              ~~~ f - \epsilon ~\leq~ T(x, D) P ~\leq~ f ~~\mbox{on}~~ \Omega \setminus
                                                             \Gamma_\epsilon
   \end{array}
\]
Furthermore, one can also assume that the Lebesgue measure of
$\Gamma_\epsilon$ is zero, namely
\[~~~ mes~ (\Gamma_\epsilon) ~=~ 0.\]
\hfill $\Box$
\end{proposition}

In view of Proposition \ref{prop}, the
spaces of piecewise smooth functions given by%
\begin{equation}\label{cndomega}
~~~ C^l_{nd} (\Omega) ~=~
      \left \{ ~~ u ~~
          \begin{array}{|l}
              \exists ~~ \Gamma \subset \Omega~~ \mbox{closed, nowhere dense}~~ : \\ \\
               ~~~~~ *)~ u : \Omega \setminus \Gamma \rightarrow {\bf R} \\ \\
               ~~~~ **)~ u \in C^l (\Omega \setminus \Gamma)
               \end{array} ~\right \}
\end{equation}
where $l \in {\bf N}$, are considered. It is easy to see that we
have the
inclusion%
\begin{equation}\label{PDEincl}
T(x, D)~ C^m_{nd} (\Omega)~\subseteq~ C^0_{nd} (\Omega)
\end{equation}

The general existence result obtained in \cite{Rosinger} is
represented through the following equation
\begin{equation}\label{PDEexists}
T(x, D)^{\#}~ ( C^m_{nd} (\Omega) )^{\#}_T ~=~ ( C^0_{nd} (\Omega)
)^{\#}
\end{equation}
Here $( C^m_{nd} (\Omega) )^{\#}_T$ and $( C^0_{nd} (\Omega)
)^{\#}$ are Dedekind order completions of $C^m_{nd} (\Omega)$ and
$C^0_{nd} (\Omega)$, respectively, when these latter two spaces
are considered with suitable partial orders. The respective
partial order on $C^m_{nd} (\Omega)$ may depend on the nonlinear
partial differential operator $T(x, D)$ in (\ref{PDEincl}), while
the partial order on $C^0_{nd} (\Omega)$ is the natural pointwise
one at the points where the two functions compared are both
continuous. The operator $T(x, D)^{\#}$ is a natural extension of
the nonlinear partial differential operator $T(x, D)$ in
(\ref{PDEincl}) to the mentioned Dedekind order completions.

Equation (\ref{PDEexists}) means that for every right hand term $f
\in ( C^0_{nd} (\Omega) )^{\#}$ in (\ref{PDE}), there exists a
solution $u \in ( C^m_{nd} (\Omega) )^{\#}_T$, and as seen later,
the set $( C^0_{nd} (\Omega) )^{\#}$ contains many discontinuous
functions beyond those piecewise discontinuous.

There is an obvious ambiguity with the piecewise smooth functions
in $C^m_{nd}(\Omega)$. Indeed, given any such function $u$, the
corresponding closed, nowhere dense set $\Gamma$ cannot be defined
uniquely. Therefore, it is convenient to factor out this
ambiguity. For the space $C^0_{nd} (\Omega)$ which is the largest
of these spaces of functions and also the range of $T(x,D)$ in
(\ref{PDEincl}) this is done by defining on it the equivalence
relation $u\sim v$ for any two elements $u,v\in C^0_{nd}
(\Omega)$, as given by
\begin{equation}
u\sim v\Longleftrightarrow\left[
\begin{tabular}
[c]{l}%
$\exists~\Gamma~\subset~\Omega$ closed, nowhere dense:\\
~~~(i) $u,v\in C(\Omega\setminus\Gamma)$\\
~~~(ii) $u=v$ $\ $on $\ \Omega\setminus\Gamma$%
\end{tabular}
\ \ \ \ \ \ \ \right]  . \label{equivrel}%
\end{equation}
The mentioned ambiguity is eliminated by going to the quotient
space
\begin{equation}\label{quspace}
\mathcal{M}^0(\Omega)=C^0_{nd}(\Omega)/\sim
\end{equation}
The partial order on $C^0_{nd}(\Omega)$ induces a partial order on
the quotient space $\mathcal{M}^0(\Omega)$, namely, for any two
$\mathbf{u},~\mathbf{v}\in\mathcal{M}^0(\Omega)$ we have
\begin{equation}
\mathbf{u}\leq\mathbf{u}\Longleftrightarrow\left[
\begin{tabular}
[c]{l}%
$\exists$ $u\in\mathbf{u},$ $v\in\mathbf{v},$
$\Gamma\subset\Omega$ closed,
nowhere dense:\\
i) $u,v\in C(\Omega\setminus\Gamma)$\\
ii) $u\leq v$ $\ $on $\ \Omega\setminus\Gamma$%
\end{tabular}
\ \ \ \ \ \ \ \right]  . \label{morder}%
\end{equation}

Using similar manipulations, this time also involving the operator
$T(x,D)$, the ambiguity in the domain of $T(x,D)$ in
(\ref{PDEincl}) is factored out, thus producing the space
$\mathcal{M}^0_T(\Omega)$ with partial order which may also depend
on the operator $T$. For details on this procedure see
\cite{Rosinger}. The equation (\ref{PDEexists})
is now replaced by%
\begin{equation}\label{PDEmexists}
T(x, D)^{\#}~ ( \mathcal{M}^m_T (\Omega) )_T^{\#} ~=~ (
\mathcal{M}^0(\Omega) )^{\#}
\end{equation}
The basic regularity result in \cite{Rosinger} is obtained by
embedding $(\mathcal{M}^0(\Omega) )^{\#}$ in the set of all
measurable functions on $\Omega$. Hence the solutions of
(\ref{PDE}) can be assimilated with measurable functions.

In the next section we will show that $\mathcal{M}^0(\Omega)$ can
be embedded in the set $\mathbb{H}(\Omega)$ of all Hausdorff
continuous functions on $\Omega$. Since the set
$\mathbb{H}(\Omega)$ is order complete it also contains the
Dedekind order completion of $\mathcal{M}^0(\Omega)$. More
precisely, we obtain that $\mathcal{M}^0(\Omega)$ is order
isomorphic to $\mathbb{H}_{nf}(\Omega)$. Hence the solutions of
(\ref{PDE}) can be assimilated with nearly finite Hausdorff
continuous functions.

\section{Assimilating the solutions of nonlinear PDEs with Hausdorff continuous
functions}

Let $u\in C_{nd}(\Omega)$. According to (\ref{cndomega}), there
exists a closed, nowhere dense set $\Gamma\subset\Omega$ such that
$u\in C(\Omega \setminus\Gamma)$. Since $\Omega\setminus\Gamma$
is open and dense in $\Omega,$ we can define%
\begin{equation}\label{F0def}
F_{0}(u)=F(\Omega\setminus\Gamma,\Omega,u)
\end{equation}
The closed, nowhere dense set $\Gamma$ used in (\ref{F0def}), is
not unique. However, we can show that the value of
$F(\Omega\setminus\Gamma,\Omega,u)$ does not depend on the set
$\Gamma$ in the sense that for every closed, nowhere dense set
$\Gamma$ such that $u\in C(\Omega\setminus\Gamma)$ the value of
$F(\Omega\setminus\Gamma,\Omega,u)$ remains the same.

Let $\Gamma_{1}$ and $\Gamma_{2}$ be closed, nowhere dense sets
such that $u\in C(\Omega\setminus\Gamma_{1})$ and $u\in
C(\Omega\setminus\Gamma_{2})$. Then the set $\Gamma_1\cup\Gamma_2$
is also closed and nowhere dense. According to Theorem
\ref{textend} in the Appendix the functions
$F(\Omega\setminus\Gamma_1,\Omega,u)$,
$F(\Omega\setminus\Gamma_2,\Omega,u)$ and
$F(\Omega\setminus(\Gamma_1\cup\Gamma_2),\Omega,u)$ are all
H-continuous and for every $x\in\Omega\setminus(\Gamma_1\cup\Gamma_2)$ we have%
\[
F(\Omega\setminus\Gamma_1,\Omega,u)(x)=
F(\Omega\setminus\Gamma_2,\Omega,u)(x)=
F(\Omega\setminus(\Gamma_1\cup\Gamma_2),\Omega,u)(x)=u(x).
\]
Since $\Omega\setminus(\Gamma_1\cup\Gamma_2)$ is dense in
$\Omega$ Theorem \ref{tindent} implies that %
\[
F(\Omega\setminus\Gamma_1,\Omega,u)=
F(\Omega\setminus\Gamma_2,\Omega,u)=
F(\Omega\setminus(\Gamma_1\cup\Gamma_2),\Omega,u).
\]
Therefore, the mapping
\[F_{0}:C_{nd}(\Omega)\longmapsto\mathbb{A}(\Omega)\]
is unambiguously defined through (\ref{F0def}). In analogy with
(\ref{gcgen}), we call $F_0$ a graph completion mapping on
$C_{nd}(\Omega)$. As mentioned above already it follows from
Theorem \ref{textend} that for every $u\in C_{nd}(\Omega)$ we have%
\[F_0(u)\in\mathbb{H}(\Omega).
\]
Furthermore, if $u\in C(\Omega\setminus\Gamma)$ we have%
\begin{equation}\label{F0inden}
F_0(u)(x)=u(x),\ x\in\Omega\setminus \Gamma.
\end{equation}
The above identity shows that the values of the function $F_0(u)$
are finite on the open and dense set $\Omega\setminus \Gamma$.
Hence, $F_0(u)$
is nearly finite, see Definition \ref{Defnearlyfinite}. Thus, we have%
\begin{equation}\label{CndHnf}
F_{0}:C_{nd}(\Omega)\longmapsto\mathbb{H}_{nf}(\Omega) %
\end{equation}

The following theorem shows that the images of two functions in
$C_{nd}(\Omega)$ under the mapping $F_0$ are the same if and only
if these functions are equivalent with respect to the relation
(\ref{equivrel}).

\begin{theorem}
\label{tfgequiv}Let $u,v\in C_{nd}(\Omega)$. Then
\[
F_0(u)=F_0(v)\ \Longleftrightarrow\ u\sim v
\]
\end{theorem}

\begin{proof}
\underline{Implication to the left.} Let $\Gamma$ be closed,
nowhere dense subset of $\Omega$ associated with $u$
and $v$ in terms of (\ref{equivrel}), that is,%
\begin{eqnarray*}
&&u,v \in C(\Omega\setminus\Gamma),\\
&&u(x)=v(x),\;x\in\Omega\setminus\Gamma.
\end{eqnarray*}
The required equality follow from (\ref{F0def}) where the set
$\Gamma$ is the one considered above. Indeed, we have
\[
F_0(u)=F(\Omega\setminus\Gamma,\Omega,u)=F(\Omega\setminus\Gamma,\Omega,v)=F_0(v)
\]

\underline{Implication to the right.} Let us denote by
$\Gamma_{1}$ and $\Gamma_{2}$ the closed, nowhere dense sets
associated with the functions $u$ and $v$, respectively, in terms
of (\ref{cndomega}), that is, $\Gamma_{1}$ and $\Gamma_{2}$ are
such that $u\in C(\Omega\setminus\Gamma_{1})$ and $v\in
C(\Omega\setminus\Gamma_{2})$. Assume that
\begin{equation}\label{assum}%
F_{0}(u)=F_{0}(v).
\end{equation}
The functions $u$ and $v$ are both continuous on the set $\Omega
\setminus(\Gamma_{1}\cup\Gamma_{2})$. Therefore, from the
property (\ref{F0inden}) and the assumption (\ref{assum}) it follows that%
\[
u(x)=F_{0}(u)(x)=F_{0}(v)(x)=v(x),{\rm  \
}x\in\Omega\setminus(\Gamma _{1}\cup\Gamma_{2}).
\]
Since the set $\Gamma_{1}\cup\Gamma_{2}$ is closed and nowhere
dense in $\Omega$, the above identity implies that $u\sim v$, see
(\ref{equivrel}).
\end{proof}

In view of (\ref{CndHnf}) and Theorem \ref{tfgequiv} now we can
define a mapping
\[
\mathbf{F}_{0}:\mathcal{M}^0(\Omega)\longmapsto\mathbb{H}_{n\!f}(\Omega)
\]
in the following way. Let $\mathbf{u}\in\mathcal{M}^0(\Omega)$ and
let $\phi
\in\mathbb{H}_{n\!f}(\Omega)$. Then%
\begin{equation}
\mathbf{F}_{0}(\mathbf{u})=\phi\Longleftrightarrow\exists u\in\mathbf{u}%
:F_{0}(u)=\phi. \label{deffg}%
\end{equation}
It is easy to see that the definition of
$\mathbf{F}_{0}(\mathbf{u})$ does not depend on the particular
representative $u$ of the equivalence class $\mathbf{u}$. Indeed,
if $u,h\in\mathbf{u}$ then $u\sim h$. Thus, $F_{0}(u)=F_{0}(h)$,
see Theorem \ref{tfgequiv}. Therefore the statement (\ref{deffg})
can be reformulated as
\begin{equation}
\mathbf{F}_{0}(\mathbf{u})=\phi\Longleftrightarrow\forall u\in\mathbf{u}%
:F_{0}(u)=\phi. \label{deffg1}%
\end{equation}

\begin{theorem}
\label{tHnfembed1}The mapping $\mathbf{F}_{0}:\mathcal{M}^0(\Omega
)\longmapsto\mathbb{H}_{n\!f}(\Omega)$ defined by (\ref{deffg}) is
an order isomorphic embedding with respect to the order relation
(\ref{morder}) in
$\mathcal{M}^0(\Omega)$ and the order relation (\ref{forder}) in $\mathbb{H}%
_{n\!f}(\Omega)$, that is, for any
$\mathbf{u}$,$\mathbf{v}\in\mathcal{M}^0(\Omega)$ we have%
\[
\mathbf{u}\leq\mathbf{v}\Longleftrightarrow\mathbf{F}_{0}(\mathbf{u}%
)\leq\mathbf{F}_{0}(\mathbf{v}).
\]
\end{theorem}

\begin{proof}
Let $\mathbf{u}$,$\mathbf{v}\in\mathcal{M}^0(\Omega)$ and $\mathbf{u}%
\leq\mathbf{v}$. According to (\ref{morder}) there exists a
closed, nowhere dense set $\Gamma$ in $\Omega$ and functions
$u\in\mathbf{u}$, $v\in \mathbf{v}$ such that $u,v\in
C(\Omega\setminus\Gamma)$ and $u(x)\leq v(x)$ for all
$x\in\Omega\setminus\Gamma.$ Using the same set $\Gamma$ in the
evaluation of $F_{0}(u)$ and $F_{0}(v)$ according to (\ref{F0def})
as well as the monotonicity of the graph completion operator, see Theorem \ref{tmon}, we have%
\[
F_0(u)=F(\Omega\setminus\Gamma,u)\leq
F(\Omega\setminus\Gamma,v)=F_0(v)
\]

Let us assume now that
$\mathbf{u},\mathbf{v}\in\mathcal{M}^0(\Omega)$ and
$\mathbf{F}_{0}(\mathbf{u})\leq\mathbf{F}_{0}(\mathbf{v})$. From
the representation (\ref{deffg}) of $\mathbf{F}_{0}$ it follows
that there exist
$u\in\mathbf{u}$ and $v\in\mathbf{v}$ such that $F_{0}(u)=\mathbf{F}%
_{0}(\mathbf{u})$ and $F_{0}(v)=\mathbf{F}_{0}(\mathbf{v}).$
Obviously we have
\begin{equation}
F_{0}(u)\leq F_{0}(v). \label{ineq9}%
\end{equation}
Since $u,v\in C_{nd}(\Omega)$ there exist closed, nowhere dense
sets $\Gamma_{1}$ and $\Gamma_{2}$ such that $u\in
C_{nd}(\Omega\setminus \Gamma_{1})$ and $v\in
C_{nd}(\Omega\setminus\Gamma_{2}).$ The set
$\Gamma=\Gamma_{1}\cup\Gamma_{2}$ is also closed, nowhere dense.
Both functions $u$ and $v$ are continuous on
$\Omega\setminus\Gamma.$ Therefore $F_{0}(u)(x)=u(x)$ and
$F_{0}(v)(x)=v(x)$ for all $x\in\Omega\setminus\Gamma $, see
(\ref{F0inden}). Hence, inequality (\ref{ineq9}) implies
\[
u(x)\leq v(x),\;x\in\Omega\setminus\Gamma,
\]
which means that $\mathbf{u}\leq\mathbf{v}$, see (\ref{morder}).
\end{proof}

\begin{theorem}
\label{tHnfembed2}Let $h\in\mathbb{H}_{n\!f}(\Omega)$. There
exists a subset
$\mathcal{G}$ of $\mathcal{M}^0(\Omega)$ such that $h=\sup\mathbf{F}%
_{0}(\mathcal{G)}$, where $\mathbf{F}_{0}(\mathcal{G)}$ is the
range of
$\mathcal{G}$ under $\mathbf{F}_{0}$, that is, $\mathbf{F}_{0}(\mathcal{G}%
)=\{\mathbf{F}_{0}(\mathbf{u)}:\mathbf{u}\in\mathcal{G}\}.$
\end{theorem}

\begin{proof}
The set $\Gamma_{n\!f}(h)=\{x\in\Omega:\infty\in h(x)$ $or$
$-\infty\in
h(x)\}$ is closed, nowhere dense, see Definition \ref{Defnearlyfinite}, and $h\in\mathbb{H}%
_{f\!t}(\Omega\setminus\Gamma_{n\!f}(h))$. Then according to
Theorem \ref{tocomplC} the function $h$ can be represented on the
set $\Omega
\setminus\Gamma_{n\!f}(h)$ as%
\begin{equation}
h(x)=(\sup\mathcal{F})(x)\ , \
x\in\Omega\setminus\Gamma_{n\!f}(h),
\label{hxsupf}%
\end{equation}
where
\[
\mathcal{F}=\{v\in C(\Omega\setminus\Gamma_{n\!f}(h)):v(x)\leq
h(x)\ , \ x\in\Omega\setminus\Gamma_{n\!f}(h)\}.
\]
The set $\mathcal{F}$ is a subset of $C_{nd}(\Omega)$ because
$\Gamma
_{n\!f}(h)$ is closed and nowhere dense. We will show that%
\[
h=\sup\mathbf{F}_{0}(\mathcal{G)}%
\]
where
\[
\mathcal{G}=\{\mathbf{v}\in\mathcal{M}^0(\Omega):\exists v\in\mathcal{F}%
:v\in\mathbf{v}\}.%
\]
Indeed, since all functions in $\mathcal{F}$\ are continuous on
$\Omega \setminus\Gamma_{n\!f}(h)$, for every
$\mathbf{v}\in\mathcal{G}$ and $v\in\mathbf{v}$ we have, see
Theorem \ref{tcont},
\begin{equation}
v(x)=F_{0}(v)(x)=\mathbf{F}_{0}(\mathbf{v})(x)\ , \ x\in\Omega
\setminus\Gamma_{n\!f}(h). \label{ident133}%
\end{equation}
Hence%
\[
\mathbf{F}_{0}(\mathbf{v})(x)=v(x)\leq h(x)\ , \
x\in\Omega\setminus \Gamma_{n\!f}(h)\}.
\]
Using that both $\mathbf{F}_{0}(\mathbf{v})$ and $h$ are
H-continuous on $\Omega$ we obtain from Theorem \ref{tindent} that
\[
\mathbf{F}_{0}(\mathbf{v})(x)\leq h(x)\ , \ x\in\Omega\ , \
\mathbf{v}\in\mathcal{G}.
\]
Therefore, $h$ is an upper bound of $\mathbf{F}_{0}(\mathcal{G)}$.
As a
bounded subset of $\mathbb{H}_{n\!f}(\Omega)$ the set $\mathbf{F}%
_{0}(\mathcal{G)}$ has a supremum in $\mathbb{H}_{n\!f}(\Omega)$,
see Theorem \ref{tocompl}. Let
$g=\sup\mathbf{F}_{0}(\mathcal{G)}$. Clearly
\begin{equation}
g\leq h. \label{ineq144}%
\end{equation}
Furthermore, from (\ref{ident133}) it follows that for every
$v\in\mathcal{F}$
and the respective class $\mathbf{v}\in\mathcal{G}$ containing $v$ we have%
\[
v(x)=\mathbf{F}_{0}(\mathbf{v})(x)\leq g(x)\ , \
x\in\Omega\setminus \Gamma_{n\!f}(h).
\]
Hence, $g$ is an upper bound of $\mathcal{F}$ on the set $\Omega
\setminus\Gamma_{n\!f}(h)$ while $h$ is the supremum of
$\mathcal{F}$ on
$\Omega\setminus\Gamma_{n\!f}(h)$. Therefore,%
\[
h(x)\leq g(x),\;x\in\Omega\setminus\Gamma_{n\!f}(h).
\]
Using the H-continuity of $g$ and $h,$ from Theorem \ref{tindent}
we obtain
that%
\[
h(x)\leq g(x),\;x\in\Omega.
\]
This together with (\ref{ineq144}) shows that $h=g=\sup\mathbf{F}%
_{0}(\mathcal{G)}$ which completes the proof.
\end{proof}

The Theorem \ref{tHnfembed2} shows that
$\mathbb{H}_{n\!f}(\Omega)$ is the smallest Dedekind order
complete subset of $\mathbb{H}(\Omega)$ which contains the image
of $\mathcal{M}^0(\Omega)$ under the order isomorphic embedding
$\mathbf{F}_{0}$. Hence it is order isomorphic to the Dedekind
order completion $\mathcal{M}^0(\Omega)^{\#}$ of
$\mathcal{M}^0(\Omega)$. The mapping discussed in this section are
illustrated on the following diagram, $\mathbf{F}_{0}{}^{\#}$
denoting the order isomorphism from $\mathcal{M}^0(\Omega)^\#$ to
$\mathbb{H}(\Omega)$.\\

\unitlength 1mm
\begin{center}
\begin{picture}(100,60)
\put(10,10){$\mathcal{M}^{0}(\Omega )^{\#}$}%
\put(90,10){$\mathbb{H}_{n\!f}(\Omega )$}%
\put(55,14){$\mathbf{F}_{0}^{\#}$}%
\put(35,11){\vector(1,0){50}}%
\put(45,6){(order isomorphism)}%
\put(15,28){\vector(0,-1){12}}%
\put(95,28){\vector(0,-1){12}}%
\put(95,23){\vector(0,1){5}}%
\put(0,31){$\mathcal{M}^{0}(\Omega )=C_{nd}^{0}(\Omega )/\sim$}%
\put(55,35){$\mathbf{F}_{0}$}%
\put(35,32){\vector(1,0){50}}%
\put(38,27){(order isomorphic embedding)}%
\put(90,31){$\mathbb{H}_{n\!f}(\Omega )$}%
\put(15,49){\vector(0,-1){12}}%
\put(10,52){$C_{nd}^{0}(\Omega )$}%
\put(90,52){$\mathbb{H}_{n\!f}(\Omega )$}%
\put(55,56){$F_0^{\#}$}%
\put(35,53){\vector(1,0){50}}%
\put(39,48){(graph completion mapping)}%
\put(95,49){\vector(0,-1){12}}%
\put(95,44){\vector(0,1){5}}%
\end{picture}
\end{center}

The set of solutions $\mathcal{M}^m_T(\Omega) )_T^{\#}$ is mapped
onto the set $\mathbb{H}_{n\!f}(\Omega )$ of all nearly finite
Hausdorff continuous functions through the composition of the
mappings $T(x,D)^{\#}$ and $\mathbf{F}_{0}^{\#}$. Considering that
both mappings are order isomorphisms, the set of solutions
$\mathcal{M}^m_T(\Omega) )_T^{\#}$ is order isomorphic with the
set $\mathbb{H}_{n\!f}(\Omega )$. Hence, the solutions of
(\ref{PDE}) through the order completion method can be assimilated
with nearly finite Hausdorff continuous functions.

\section{Conclusion}
The paper deals with the regularity of the solutions of nonlinear
PDEs obtained through the order completion method. We show that
these solutions can be assimilated with Hausdorff continuous
function, thus significantly improving the results in
\cite{Rosinger} with respect to the regularity properties of the
solutions. The applications of the class of Haudorff continuous
functions discussed here as well as in other recent publications,
\cite{QM}, \cite{Intro}, \cite{Anguelov-Rosinger}, show that this
class may play an important role in what is typically called Real
Analysis. In particular, one may note that one of the main engines
behind the development of the various spaces in Real and Abstract
Analysis are the partial differential equations with the need to
assimilate the various types of ''weak'' solutions. Since the
solutions of very large classes of nonlinear partial differential
equations can be assimilated with nearly finite Hausdorff
continuous functions, the set of these functions might be a viable
alternative to some of the presently used functional spaces (e.g.
$L^{p}(\Omega),$ Sobolev spaces) with the advantage of being both
more regular and universal.

\section*{Appendix}
The concept of Hausdorff continuous interval valued functions was
developed first within the theory of Hausdorff approximations of
real functions, see \cite{Sendov}. The name is derived from the
fact that for a Hausdorff continuous function
$f=[\underline{f},\overline{f}]$ the Hausdorff distance between
the graphs of $\underline{f}$ and $\overline{f}$ is zero. Since
the Hausdorff continuous functions are in general interval valued
they are also studied as a part of the Interval Analysis, see
\cite{Anguelov-Markov}, \cite{Intro}.

The minimality condition associated with the Hausdorff continuity,
see Definition \ref{DefHcont}, requires that the graph of a
Hausdorff continuous function is as 'thin' as possible, that is,
the function assumes proper interval values only when it is
necessary to ensure that the graph of this interval function is a
closed subset of $\Omega\times\overline{\mathbb{R}}$. As a result
the set where a Hausdorff continuous function assumes proper
interval values is small. The next theorem shows that this set is
meager or a set of first Baire category, that is, a countable
union of closed and nowhere dense sets.

\begin{theorem}
\label{FBC}The set $W_f=\{x\in \Omega:w(f(x))>0\}$ of all points
where $f\in\mathbb{A}(\Omega)$ assumes proper interval values is a
set of first Baire category.
\end{theorem}

It may appear at first that the minimality condition in Definition
\ref{DefHcont} applies at each individual point $x$ of $\Omega$,
thus, not involving neighborhoods. However, the graph completion
operator $F$ does appear in this condition. And this operator
according to (\ref{gcgen}) and therefore (\ref{lbfgen}) and
(\ref{ubfgen}) does certainly refer to neighborhoods of points in
$\Omega$, a situation typical, among others, for the concept of
continuity. Hence the following property of the continuous
functions is preserved.

\begin{theorem}
\label{tindent}Let $f,g$ be H-continuous on $\Omega$ and let $D$
be a dense subset
of $\Omega$. Then%
\begin{eqnarray*}
&{\rm a)}& \ \ f(x)~\leq~ g(x),\;x\in D~\Longrightarrow~
f(x)~\leq~
g(x),\;x\in \Omega,{\rm  \ \ \ \ \ \ \ \ \ \ \ \ \ \ \ \ \ \ \ \
}\\
&{\rm b)}& \ \ f(x)~=~g(x),\;x\in D~\Longrightarrow~ f(x)~=~ g(x),\;x\in \Omega.{\rm  \ \ \ \ \ \ \ \ \ \ \ \ \ \ \ \ \ \ \ \ }%
\end{eqnarray*}
\end{theorem}

The following two theorems represent essential links with the
usual point valued continuous functions.

\begin{theorem}\label{tcont}
Let $f=[\underline{f},\overline{f}]$ be an
H-continuous function on $\Omega$. \newline a) If $\underline{f}$
or $\overline{f}$ is continuous at a point $a\in \Omega$ then
$\underline{f}(a)=\overline{f}(a)$.\newline b) If
$\underline{f}(a)=\overline{f}(a)$ for some $a\in \Omega$ then
both $\underline{f}$ and $\overline{f}$ are continuous at $a$.
\end{theorem}

\begin{theorem}\label{textend}
Let $D$ be a dense subset of $\Omega$. If $f\in C(D)$ then
\begin{eqnarray*}
&&F(D,\Omega,f)\in\mathbb{H}(\Omega),\\
&&F(D,\Omega,f)(x)=f(x),\ x\in D.
 \end{eqnarray*}
\end{theorem}

A partial order which extends the total order on
$\overline{\mathbb{R}}$ can be defined on
$\mathbb{I\,}\overline{\mathbb{R}}$ in more than one way.
Historically, several partial orders are associated with the set
$\mathbb{I\,}\overline{\mathbb{R}}$, namely, \newline $(i)$ \ the
inclusion
relation $[\underline{a},\overline{a}]\subseteq\lbrack\underline{b}%
,\overline{b}]\Longleftrightarrow\underline{b}\leq\underline{a}\leq
\overline{a}\leq\overline{b}$\newline $(ii)$ \ the ''strong''
partial order
$[\underline{a},\overline{a}]\preceq\lbrack\underline{b},\overline
{b}]\Longleftrightarrow\overline{a}\leq\underline{b}$\newline
$(iii)$ \ the partial order defined by (\ref{iorder}).\newline The
use of the inclusion relation on the set
$\mathbb{I\,}\overline{\mathbb{R}}$ is motivated by the
applications of interval analysis to generating enclosures of
solution sets. However, the role of partial orders extending the
total order on $\overline{\mathbb{R}}$ has also been recognized in
computing, see \cite{Birkhoff}. Both orders $(ii)$ and $(iii)$ are
extensions of the order on $\overline{\mathbb{R}}$. The use of the
order $(ii)$ is based on the view point that inequality between
intervals should imply inequality between their interiors. This
approach is rather limiting since the order $(ii)$ does not retain
some essential properties of the order on $\overline{\mathbb{R}}$.
For instance, a proper interval $A$ and the interval
$A+\varepsilon$ are not comparable with respect to the order
$(ii)$ when the positive real number $\varepsilon$ is small
enough. The partial order $(iii)$ is introduced and studied by
Markov, see \cite{Markov}, \cite{Markov2}. The results reported in
\cite{QM} and in the present paper indicate that indeed the
partial order (\ref{forder}) induced pointwise by (\ref{iorder})
is an appropriate partial order to be associated with the
Hausdorff continuous interval valued functions.

The monotonicity with respect to the relation inclusion was
discussed in Section 2, see (\ref{Finclmon1}) and
(\ref{Finclmon2}). The following theorem states the monotonicity
of the Baire operators and the graph completion operator with
respect to the order (\ref{forder}) induces in a pontwise way by
the order (\ref{iorder}).

\begin{theorem}\label{tmon}
The lower Baire operator, the upper Baire operator and the graph
completion operator are all monotone increasing with respect to
the order (\ref{forder}) on the respective domains and ranges,
that is, if $D$ is a dense subset of
$\Omega$, for every two functions $f,g\in\mathbb{A}(D)$ we have%
\[
f(x)\leq g(x),~x\in D~\Longrightarrow ~ \left\{\begin{tabular}{l}
$I(D,\Omega,f)(x)\leq I(D,\Omega,g)(x),~x\in\Omega$\\
$S(D,\Omega,f)(x)\leq S(D,\Omega,g)(x),~x\in\Omega$\\
$F(D,\Omega,f)(x)\leq F(D,\Omega,g)(x),~x\in\Omega$
\end{tabular}\right.
\]
\end{theorem}

Important property of the set $\mathbb{H}(\Omega)$ is that it is
order complete. As we noted, the order completeness or the
Dedekind order completeness is not a typical property for the
spaces of functions considered in Real Analysis. In this way, the
class of H-continuous functions and its subclasses mentioned in
the next theorem can provide solutions to open problems or improve
earlier results related to order.

\begin{theorem}.\label{tocompl}

\begin{itemize}
\item[(i)]The set $\mathbb{H}(\Omega)$ of all H-continuous
functions is order complete.

\item[(ii)]The set $\mathbb{H}_{bd}(\Omega)$ of all bounded
H-continuous
functions, that is,%
\[
\mathbb{H}_{bd}(\Omega)=\{f\in\mathbb{H}(\Omega):\exists
M\in\mathbb{R}:|f(x)|\leq M,\ x\in\Omega\}
\]
is Dedekind order complete

\item[(iii)]The set $\mathbb{H}_{f\!t}(\Omega)$ of all finite
H-continuous functions, that is,
\[
\mathbb{H}_{f\!t}(\Omega)=\{f\in\mathbb{H}(\Omega):|f(x)|<+\infty,\
x\in\Omega\}
\] is Dedekind order complete

\item[(iv)]The set $\mathbb{H}_{n\!f}(\Omega)$ of all nearly
finite H-continuous functions, that is,
\[
\mathbb{H}_{nf}(\Omega)=\{f\in\mathbb{H}(\Omega):\exists D-\ {\rm
dense\ subset\ of}\ \Omega:|f(x)|<+\infty,\ x\in D\}
\]
is Dedekind order complete.
\end{itemize}
\end{theorem}

The resent paper \cite{QM} gives the Dedekind order completion of
the space $C(X)$ of all continuous real functions on a topological
space $X$ in terms of Hausdorff continuous functions, thus
improving significantly an earlier result by Dilworth, see
\cite{Dilworth}. The main result in \cite{QM} is stated below for
the case when $X=\Omega$.

\begin{theorem}\label{tocomplC}
The set $\mathbb{H}_{f\!t}(\Omega)$ is a Dedekind order completion
of the set $C(\Omega)$. Moreover, for every
$h\in\mathbb{H}_{f\!t}(\Omega)$ we have
\[
h=sup\{f\in C(\Omega):f\leq h\}
\]
\end{theorem}

The proof of the theorems in this apendix can be found in
\cite{QM} and \cite{Intro}.

\end{document}